\documentclass[12pt,centertags,draft]{amsart}

\usepackage{enumerate}
\usepackage{amssymb}
\usepackage{amsfonts}
\usepackage{amsthm}
\usepackage{amsmath}
\usepackage[mathscr]{eucal}
\usepackage{latexsym}

\numberwithin{equation}{section}

\DeclareMathOperator*{\slim}{s-lim}
\DeclareMathOperator*{\nlim}{n-lim}

{\bf}{\it}
{\bf}{\it}

\theoremstyle{plain}
\newtheorem{thm}{Theorem}[section]

\newtheorem{lemma}[thm]{Lemma}

\newtheorem{cor}[thm]{Corollary}
\newtheorem{hyp}[thm]{Hypothesis}
\newtheorem{rem}[thm]{Remark}

\theoremstyle{definition}
\newtheorem{defn}[thm]{Definition}

\theoremstyle{remark}

\newcommand{\la}{\lambda}

\newcommand{\De}{\Delta}

\renewcommand{\epsilon}{\varepsilon}

\newcommand{\R}{\mathbb{R}}

\newcommand{\Z}{\mathbb{Z}}

\newcommand{\C}{\mathbb{C}}

\newcommand{\Hh}{\mathcal{H}}
\newcommand{\he}{\mathcal{H}}
\newcommand{\K}{\mathcal{K}}
\newcommand{\D}{\mathcal{D}}
\newcommand{\A}{\mathcal{A}}

\newcommand{\cH}{{\mathcal H}}

\renewcommand{\Re}{\text{\rm Re}}
\renewcommand{\Im}{\text{\rm Im}}

\newcommand{\im}{\text{\rm Im}}
\renewcommand{\i}{\text{\rm i}}
\newcommand{\ii}{\text{\rm i}}

\newcommand{\Dim}{\text{\rm Dim}}

\renewcommand{\index}{\text{\rm ind}}

\newcommand{\fM}{\mathfrak{M}}

\newcommand{\Iinfty}{\mathrm{I_\infty}}

\newcommand{\IIone}{\mathrm{II_1}}

\title[The Birman-Schwinger principle]{The Birman-Schwinger principle in
von Neumann algebras of finite type}

\author[V.~Kostrykin]{Vadim Kostrykin}
\address{Vadim Kostrykin\\ Institut f\"{u}r
Mathematik, Technische Universit\"{a}t Claus\-thal, Erzstra{\ss}e 1,
D-38678 Clausthal-Zellerfeld, Germany}
\email{kostrykin@math.tu-clausthal.de,\hskip -0.005cm kostrykin@t-online.de}

\author[K.~A.~Makarov]{Konstantin A.~Makarov}
\address{Konstantin A.~Makarov\\ Department of Mathematics\\ University of Missouri\\
Co\-lum\-bia, MO 65211, USA} \email{makarov@math.missouri.edu}

\author[A.~Skripka]{Anna Skripka}
\address{Anna Skripka\\ Department of Mathematics\\ University of Missouri\\
Co\-lum\-bia, MO 65211, USA} \email{skripkaa@math.missouri.edu}



\begin{document}

\begin{abstract}
We introduce a relative index for a pair of dissipative operators in
a von Neumann algebra of finite type and prove an analog of the
Birman-Schwinger principle in this setting. As an application of this
result, revisiting the Birman-Krein formula in the abstract
scattering theory, we represent the de la Harpe-Skandalis determinant
of the characteristic function of dissipative operators in the
algebra in terms of the relative index.
\end{abstract}

\maketitle


\section{Introduction }\label{sec:intro}

In 1961 M.~Sh.~Birman \cite{BirmanSw} and J.~Schwinger
\cite{Schwinger} independently introduced a method to control the
number of negative eigenvalues of Schr\"{o}dinger operators. In the
abstract operator-theoretic setting, the classical Birman-Schwinger
principle (in its simplest form) states (see, e.g., \cite{Simon}):
{\it Given a self-adjoint  strictly positive operator $H_0$ and a
non-negative self-adjoint compact operator $V$ on a Hilbert space
$\mathcal{H}$, the number of negative eigenvalues $($counting
multiplicity$)$ of the operator $H=H_0-V$ coincides with the number of
eigenvalues greater than one of the Birman-Schwinger operator
$V^{1/2}H_0^{-1}V^{1/2}$}. That is,
\begin{equation}\label{bss}
\dim[E_{H}(\R_-)\mathcal{H}]=\dim
\left[E_{I-V^{1/2}H_0^{-1}V^{1/2}}(\R_-)\mathcal{H}\right],
\end{equation}
where $\R_-=(-\infty,0)$ and $E_T(\cdot)$ is the spectral measure of a
self-adjoint operator $T$.

The sign-definiteness assumptions upon $H_0$ and $V$ can be relaxed
and the principle admits further generalizations. Assume that  $V$ is
factorized in the form $V=K^* N^{-1} K$, with $N$ a self-adjoint
unitary and $K$ a compact operator,  and  that $H_0$ and $H=H_0 - V$
have  bounded inverses. Then \eqref{bss} can be extended to a more
general equality
\begin{equation}\label{iniintr}
\index\left(E_{H_0}(\R_-),E_{H}(\R_-)\right)=
\index\left(E_{N}(\R_-),E_{N-KH_0^{-1}K^*}(\R_-)\right)
\end{equation}
of  the Fredholm indices for the associated  pairs of the spectral
projections (cf.~\cite{MakarovP} for the proof in the case of trace class
perturbations; see also \cite{Push1}, \cite{Push2}). For the concept of the
Fredholm index for a pair of orthogonal projections we refer to
\cite{Avron}.

The main purpose of this paper is to find an appropriate
generalization of the principle \eqref{iniintr} in the context of
perturbation theory in a von Neumann algebra $\A$ of finite type. To
accomplish this goal, we introduce the concept of a \emph{relative}
index $\xi(M,N)$ associated with a pair $(M,N)$ of dissipative
elements in $\A$ via
\begin{equation}\label{xi:def}
\xi(M,N) = \tau[\Xi(N)] - \tau[\Xi(M)].
\end{equation}
Here $\Xi(M)$ denotes the $\Xi$-operator \cite{MakarovL},
\cite{MakarovP} (cf. also \cite{Carey:0}) associated with $M$ and
$\tau$ a normal tracial state on the algebra $\A$.

If both $M$ and $N$ are self-adjoint, the relative index $\xi(M,N)$
can be expressed in terms of the $\tau$-Fredholm indices of the
corresponding spectral projections:
\begin{equation}\label{xi:sa}
\xi(M, N)=\index_\tau\big(E_{N}(\R_-),E_M(\R_-)\big)+\frac{1}{2}
\index_\tau\big(E_{N},E_M\big).
\end{equation}
Recall that the notion of the $\tau$-Fredholm index for a pair of
orthogonal projections $(P,Q)$ is an analog of the index introduced in
\cite{Avron}, where the usual trace has to be replaced by the tracail
state  $\tau$. In the particular case of von Neumann algebras of
finite type, one has $\index_\tau(P,Q) = \tau(P-Q)$. We refer to
\cite{Breuer1}, \cite{Breuer2} for the theory of $\tau$-Fredholm
operators.

The main result of the present paper (see Theorem \ref{mbs})
establishes a generalization of the Birman-Schwinger principle to the
case of dissipative operators in a finite  von Neumann algebra $\A$.
For boundedly invertible dissipative operators $M$, $N$,
$M-K^*N^{-1}K$, and $N-KM^{-1}K^*$ in $\A$ we prove the relation
\begin{equation}\label{BSP}
\xi(M,M-K^*N^{-1}K)= \xi(N,N-KM^{-1}K^*).
\end{equation}
In the self-adjoint
case, this relation together with \eqref{xi:sa} provides an analog of
\eqref{iniintr} for the  $\tau$-Fredholm indices.

Relaxing the invertibility assumption on the operators $M$ and/or
$N$, we present an extension of the principle \eqref{BSP} (see
Theorem \ref{bsp}). In particular, if $N$ has a bounded inverse and
the family of the operators  $N-K(M+\ii\epsilon I)^{-1}K^*$,
$\epsilon > 0$, has a limit as $\epsilon\downarrow 0$ in the norm
topology as an invertible (dissipative) operator we show that  the
relation
\begin{equation}\label{vonvon}
\xi(M, M-K^*N^{-1}K) = \xi\left(N, N-K(M+\i0 I)^{-1}K^*\right)
\end{equation}
holds.

As an application of  \eqref{vonvon} to  the self-adjoint case, we
study the perturbation problem $H_0\mapsto H=H_0-K^*N^{-1}K$, with
$N=N^*$ boundedly invertible and $H_0=H_0^*$. Under mild additional
assumptions, \eqref{vonvon} leads to an equality (see Theorem
\ref{bkf}) relating the index $\xi(H,H_0)$ to the de la
Harpe-Skandalis determinant \cite{Skandalis} of the Lifshits
characteristic function of the dissipative operator $N - K (H_0 + \ii
0 I)^{-1} K^*$. We remark that this result strongly resembles the
Birman-Krein formula \cite{BirmanK} relating the scattering matrix to
the spectral shift function.

It should be mentioned that in the context of perturbation theory for
self-adjoint operators in von Neumann algebras of finite type, the
function
\begin{equation}\label{ssf}
\R\ni\lambda\mapsto \xi(H-\lambda I, H_0-\lambda I)
\end{equation}
coincides with the spectral shift function associated with the pair
of self-adjoint operators ($H$, $H_0$). We recall that the concept of
the spectral shift function was introduced by I.~M.~Lifshits
\cite{Lifshits} and M.~G.~Krein \cite{Krein1} for (finite or
infinite) factors of type I (see \cite{BirmanP}, \cite{BirmanY},
\cite{Yafaev}, and references therein) and it  has been extended to
the case of (semi)finite von Neumann algebras in \cite{Azamov} and
\cite{Carey} (see also \cite{Boyadzhiev}).

Throughout the paper we assume that $\A$ is a von Neumann algebra of
finite type and $\tau$ a normal tracial state on it. In the
case when $\A$ is a factor of type $\IIone$, the symbol
$\text{Dim}(\cdot)$ stands for the relative dimension associated with
$\A$. The set of the boundedly invertible dissipative operators in
$\A$ is given particular  consideration and we reserve
 the symbol $\D_\A$ for this set. We use the letter $K$ to refer to an
arbitrary operator in $\A$ and $M$, $N$ to refer to dissipative
operators in $\A$. We denote self-adjoint operators in $\A$ by $H_0$,
$V$, and $H$ while  discussing issues of perturbation problems.
Auxiliary self-adjoint operators will be denoted by $A$, $B$, $L$ and
unitary operators by $U$, $S$.

\section{The $\Xi$-operator}\label{s2}

Suppose  $M$ is a dissipative, not necessarily invertible, operator
in $\A$ and $L$ its  minimal self-adjoint dilation (see
\cite{Nagy:Foias}) in a Hilbert space $\mathcal{K}\supset \Hh$. We
define the $\Xi$-operator associated with $M$ by
\begin{equation}
\Xi(M)=P_\Hh\left[E_L(\R_-)+\frac12E_L(\{0\})\right]\bigg|_{\Hh},
\end{equation}
where $E_L(\cdot)$ stands for  the spectral measure of $L$ and
$P_{\Hh}$ for  the orthogonal projection in the space $\mathcal{K}$
onto $\Hh$.

\begin{thm}\label{xiop}
If $M\in \A$ is a dissipative operator, then the self-adjoint
non-negative contraction $\Xi(M)$ belongs to the algebra $\A$.
\end{thm}

\begin{proof}
Suppose first that  $M$ has a bounded inverse. Then, by  the Langer
lemma, the minimal self-adjoint dilation $L$ of $M$ has a trivial
kernel, that is,
\begin{equation*}
E_L(\{0\})=0,
\end{equation*}
and, therefore,
\begin{equation}\label{Xi:inv}
\Xi(M)=P_\Hh E_L(\R_-)\big|_{\Hh}=\frac1\pi \Im \log M
\end{equation}
(cf.~Lemma 2.7 in \cite{MakarovL}). Here $\log M$  denotes the
principal branch of the operator logarithm of $M\in\D_\A$ with the
cut along the negative imaginary semi-axis provided by the Riesz
functional calculus. Equivalently, the operator logarithm $\log M$
can be understood as the norm-convergent Riemann integral
\begin{equation}\label{intrep}
\log M=-\i\int_0^\infty\left((M+\i\la I)
^{-1}-(1+\i\la)^{-1}I\right)d\la .
\end{equation}
Representation \eqref{intrep} proves that  $\Xi(M)$ is an element of
$\A$ (under the assumption that $M$ has a bounded inverse).

To prove the claim of the theorem in the general case, it suffices to
deduce that
\begin{equation}\label{SOT}
\slim_{\varepsilon\downarrow 0}\Xi(M+\i\epsilon I)=\Xi(M)
\end{equation}
whenever $M$ is dissipative. Indeed, given $\epsilon>0$, the
dissipative operator $M+\i\epsilon I\in \A$ obviously has  a bounded
inverse. Hence $\Xi(M+\epsilon\ii I)\in \A$, by the first part of the
proof, and the claim follows from \eqref{SOT}.

In order to prove \eqref{SOT}, we note that
\begin{equation}\label{al11}
\begin{split} & \Xi(M+\i\epsilon I) =\frac1\pi\Im\log(M+\i\varepsilon I)\\  & \qquad =
-\frac1\pi\int_0^\infty\Re((M+\i\varepsilon I+\i\la I)^{-1}-(1+\i\la)^{-1}I_{\Hh})d\la\\
 &\qquad =-\frac1\pi P_{\Hh}\int_0^\infty\Re((L+\i\varepsilon I+\i\la
I)^{-1}-(1+\i\la)^{-1}I_ {\mathcal{K}})\big|_{\Hh}d\la\\ &\qquad
=\frac1\pi P_{\Hh}\Im\log(L+\i\varepsilon I)\big|_{\Hh}
\\  &\qquad =\frac1\pi P_{\Hh}\Im\log(L+\i\varepsilon
I)\big[E_L(\{0\})+E_L(\R\setminus\{0\})\big]\big|_{\Hh}.
\end{split}
\end{equation}

Following almost verbatim the arguments in \cite{MakarovL}, we verify
that
\begin{equation}\label{al22}
\slim_{\varepsilon\downarrow 0}\frac1\pi \Im \log(L+\i\varepsilon
I)E_L(\R\setminus\{0\})=E_L(\R_-),
\end{equation} with the limit
taken in the strong operator topology of the Hilbert space $\K$.
Finally, applying the Spectral Theorem to the self-adjoint operator
$L$ with the use of \eqref{al11} and \eqref{al22}, we conclude that
\begin{align*}
\slim_{\varepsilon\downarrow 0}\Xi (M+\i\varepsilon
I)=P_{\Hh}\left[E_L(\R_-)+\frac12
E_L(\{0\})\right]\bigg|_{\Hh}=\Xi(M).
\end{align*}
\end{proof}

\begin{rem}\label{nepr}
As one can see from the proof of Theorem \ref{xiop}, the
$\Xi$-operator possesses the continuity property in the sense that
\begin{equation}
\slim_{\varepsilon\downarrow 0}\Xi(M+\i\epsilon I)=\Xi(M),
\end{equation}
whenever $M$ is a dissipative operator in $\A$. It is also  clear
that if, in addition, $M$ is self-adjoint, then the $\Xi$-operator
can be expressed in terms of the spectral resolution $E_M(\cdot)$
associated with $M$ via
\begin{equation}
\Xi(M)=E_M(\R_-)+\frac12E_M(\{0\}).
\end{equation}
\end{rem}

To conclude this section, we link the trace of the $\Xi$-operator to
the phase of the de la Harpe-Skandalis determinant \cite{Skandalis}.
Basic properties of  this determinant can be found in Appendix
\ref{sa}.

\begin{thm}
Assume that $M\in\D_\A$. Let ${\det}_{\tau}M$ be the de la
Harpe-Skandalis determinant associated with the homotopy class of the
$C^1$-paths of invertible operators joining $M$ with  the identity
$I$ and containing any $C^1$-path  $[0,1]\ni t\mapsto M_t\in \D_\A$.
Then
\[{\det}_{\tau}M=\exp(\i\pi\tau [\Xi(M)])\cdot \De(M),\]
 with $\De(\cdot)$
the Fuglede-Kadison determinant $($cf.~\cite{Kadison}$)$.
\end{thm}

\begin{proof}
As any complex number, ${\det}_{\tau}M$ can be written in the polar
form
\begin{equation}\label{polar}{\det}_{\tau}M=\exp(\i\im\log
[{\det}_{\tau}M])\cdot |{\det}_{\tau}M|.\end{equation} Lemma \ref{a2} (i) implies that
${\det}_{\tau}M=\exp(\tau[\log M])$ and Lemma \ref{a1} (ii) that
$|{\det}_{\tau}M|=\De(M)$. Combining the latter representations with
\eqref{polar}, one gets
\begin{equation}\label{ddeett}{\det}_{\tau}M=\exp(\i\pi\im\,\tau[\log M])\cdot \De(M).
\end{equation}
By positivity of the state $\tau$, one concludes that  $\tau\circ \im
=\im \circ \tau$, and hence the right hand side of \eqref{ddeett}
equals $\exp(\i\pi\tau[\im\log M])\cdot \De(M)$.  Taking into account
\eqref{Xi:inv} completes the proof.
\end{proof}

\begin{rem}\label{mors}
In the case $\A$ is a finite type factor,
$$
\tau[\Xi(H)]=\Dim [ E_H((-\infty, 0)) \he]
$$
whenever $H$ is a self-adjoint invertible element in $\A$. Thus,
$\tau [\Xi(M)]$ can be considered  a natural generalization of the
Morse index of the dissipative element $M$.
\end{rem}

\section{The Birman-Schwinger Principle}\label{s3}

The main aim of this section is to provide an analog of the
Birman-Schwinger principle in the context of perturbation theory for
dissipative operators in the von Neumann algebra setting.

\begin{defn}
We define the $\xi$-index associated with the pair $(M,N)$ of
dissipative operators $M$ and $N$  in $\A$ by
\begin{equation}\label{xidef}
\xi(M,N)=\tau[\Xi(N)]-\tau[\Xi(M)].
\end{equation}
\end{defn}

\begin{rem}
The index $\xi(M,N)$ can also be recognized as the argument of the de
la Harpe-Skandalis determinant  $\De(t\mapsto M_t)$ associated with
the homotopy class of the nonsingular $C^1$-paths joining $M$ with
$N$ and containing any $C^1$-path $[0,1]\ni t\mapsto M_t\in\D_\A$
with the endpoints  $M_0=M$ and $M_1=N$.  That is,
\[\De(t\mapsto M_t)=\exp\big ( \i \pi \xi(M,N)\big )\cdot \Delta\big (NM^{-1}\big ),\]
with $ \Delta(\cdot)$ the Fuglede-Kadison determinant.
\end{rem}

We note that in view of Remark \ref{nepr} the relative index
associated with the pair $(H_0,H)$ of self-adjoint operators in $\A$
admits a transparent representation via the $\tau$-Fredholm indices
of the corresponding spectral projections
\begin{equation}\label{xi:xi}
\xi(H,H_0)=\index_\tau\big(E_{H_0}(\R_-),E_H(\R_-)\big)+\frac{1}{2}
\index_\tau\big(E_{H_0}(\{ 0\}),E_H(\{0\})\big).
\end{equation}

We start with an invariance principle for the $\xi$-index associated
with a pair of boundedly invertible  dissipative operators, a natural
analog  of the {\it  Birman-Schwinger principle} in the perturbation theory
for self-adjoint operators in the standard $\Iinfty$ setting.

\begin{thm}\label{mbs}
Let $K\in \A$ and $M$, $N\in \D_\A$. Suppose, in addition, that the
dissipative operators $M-K^*N^{-1}K$ and $N-KM^{-1}K^*$ are boundedly
invertible. Then
\begin{equation}
\label{mbsformula} \xi(M,M-K^*N^{-1}K)= \xi(N,N-KM^{-1}K^*).
\end{equation}
\end{thm}

Before turning to the proof of Theorem \ref{mbs}, let us interpret
its result in the context of perturbation theory for self-adjoint
operators.

Assume that $H_0=H_0^*$ is a boundedly invertible element in $\A$,
and that the perturbation $V=H-H_0$ can be factored in the
form\footnote{Such a factorization is available for any $V=V^*$; for
instance, one can take $K=\sqrt{|V|}$ and $N=-\text{sgn}(V)$,
with $\text{sgn}(x)=\begin{cases}\;\;\,1 & \text{if}\quad x\geq 0,\\
-1 &\text{if}\quad x<0.\end{cases}$} $V=-K^*N^{-1}K$, with $N=N^*$ a
boundedly invertible element in $\A$. Then Theorem \ref{mbs}
guarantees the coincidence of the $\tau$-Fredholm indices for the
dual pairs of the spectral projections
\begin{equation}\label{ini}
\index_\tau\left(E_{H_0}(\R_-),E_{H}(\R_-)\right)=
\index_\tau\left(E_{N}(\R_-),E_{N-KH_0^{-1}K^*}(\R_-)\right).
\end{equation}

In particular, if $\A$ is a factor of finite type, $H_0$ and $V$ are
positive, and both $H_0$ and $H$ have bounded inverses, principle
\eqref{ini} acquires the traditional ``counting dimensions" flavor
(cf.~\eqref{bss}):
\begin{equation}\label{bspg}
\text{\rm Dim}\big[E_{H_0-V}(\R_-)\he\big]=\text{\rm
Dim}\left[E_{V^{1/2}H_0^{-1}V^{1/2}}((1,\infty))\he\right].\end{equation}

\begin{proof}[Proof of Theorem \ref{mbs}]
Introduce an auxiliary Herglotz operator-valued function
\begin{equation*}
z\mapsto{\fM}(z)=
\begin{pmatrix} M+zI& K^*\\
K& N+zI
\end{pmatrix},\quad z\in\C_+,
\end{equation*} with values in the von Neumann algebra $\A\overline{\otimes}M_2$,
where $M_2$ is the space of
$2\times 2$ (scalar) matrices. Note that ${\fM}(z)$, $z\in\C_+$,
are boundedly invertible operators in $\A\overline{\otimes}M_2$ and
the diagonal entries of ${\fM}^{-1}(z)$ are the inverses of the
operators
\begin{align}\label{shur}
\mathcal{M}(z) &= M+zI-K^*(N+zI)^{-1}K,\\
\nonumber  \mathcal{N}(z) &
=N+zI-K(M+zI)^{-1}K^*,
\end{align}
the Schur complements of ${\fM}(z)$.

Taking into account that
\[\frac{d}{dz}{\fM}(z)=\begin{pmatrix}
I& 0\\0& I\end{pmatrix}\] and using the Dixmier-Fuglede-Kadison
differentiation formula (cf. \cite{Dixmier} and \cite{Kadison})
yields
\begin{equation}
\label{3.5} \frac{d}{dz} \tau^{(2)}\big[\log
{\fM}(z)\big]=\tau^{(2)}\big[{\fM}^{-1}(z)\big]=\frac12\tau\big[{\mathcal
M}^{-1}(z)\big]+\frac12\tau\big[{\mathcal N}^{-1}(z)\big].
\end{equation}
Here $ \tau^{(2)}$ denotes the normal tracial state on the von
Neumann algebra $\A\overline{\otimes}M_2$ given by
\begin{equation}
\tau^{(2)}\left[\begin{pmatrix}
A & B\\
C & D
\end{pmatrix}\right
]= \frac{\tau(A)+\tau(D)}{2}, \quad A, B, C, D\in \A.
\end{equation}

By direct computations, we get
\begin{equation}
\label{3.6} {\mathcal N}^{-1}(z)=(N+zI)^{-1}+(N+zI)^{-1}K{\mathcal
M}^{-1}(z)K^*(N+zI)^{-1}.
\end{equation}
Employing the additivity and cyclicity of the  state $\tau$
and representation \eqref{3.6}, we derive
\begin{align}\label{3.7}
\nonumber &\tau\big[{\mathcal M}^{-1}(z)\big]+\tau\big[{\mathcal
N}^{-1}(z)\big]\\\nonumber&=\tau\big[{\mathcal
M}^{-1}(z)\big]+\tau\big[{\mathcal
M}^{-1}(z)K^*(N+zI)^{-2}K\big]+\tau\big[(N+zI)^{-1}\big]\\
\nonumber&=\tau\big[{\mathcal M}^{-1}(z)(I+K^*(N+zI)^{-2}K)\big]
+\tau\big[(N+zI)^{-1}\big]\\
&=\frac{d}{dz}\big (\tau[\log \mathcal{M}(z)]+\tau[\log (N+zI)]\big
).
\end{align}
 Comparing  \eqref{3.5} and \eqref{3.7} gives
\begin{equation*}
\frac{d}{dz}\big(2\tau^{(2)}[\log{\fM}(z)]\big)=\frac{d}{dz}\big(\tau[\log{\mathcal
M}(z)]+\tau[\log(N+zI)]\big).
\end{equation*}
{}From this, we conclude that
\begin{equation}\label{3.10}
2\tau^{(2)}[\log{\fM}(z)]=\tau[\log{\mathcal M}(z)]+\tau[\log
(N+zI)]+C,\end{equation} with $C$ a constant. Combining the
asymptotic expansions
\begin{equation*}
\begin{split} \tau^{(2)}[\log{\fM}(\ii y)] &= \log (\ii y)
+\mathcal{O}\left(\frac1y\right),\\
\tau[\log{\mathcal M}(\ii y)] &= \log (\ii
y)+\mathcal{O}\left(\frac1y\right),\\
\tau[\log(N+\ii yI)] &= \log (\ii y)+\mathcal{O}\left(\frac1y\right)
\end{split}
\end{equation*}
as $y\to +\infty$, we infer that the constant $C$ in \eqref{3.10}
equals zero and, hence,
\begin{equation}\label{kkk}
2\tau^{(2)}[\log{\fM}(z)]=\tau[\log{\mathcal M}(z)]+\tau[\log
(N+zI)], \quad \Im z>0.
\end{equation}

Computing the normal boundary values as $z\downarrow 0$ in
\eqref{kkk} ensures the equality
\begin{align}\label{3.333}
\tau^{(2)}[\log{\bf M}]=
\frac12\left(\tau[\log(M-K^*N^{-1}K)]+\tau[\log N]\right),
\end{align}
where
$$
{\bf M}=\fM(0).
$$

Next, we note that
\begin{equation*}
\begin{pmatrix} 0& I\\I& 0\end{pmatrix}
\begin{pmatrix} M& K^*\\K& N\end{pmatrix}
\begin{pmatrix} 0& I\\I& 0\end{pmatrix}^{-1}=\begin{pmatrix} N&
K\\K^*& M\end{pmatrix}.
\end{equation*}
It is straightforward to check that, for any unitary operator $U\in
\A\overline{\otimes}M_2$ and any $H\in\D_{\A\overline{\otimes}M_2}$,
\begin{equation}\label{1.7}
\log\left(UHU^{-1}\right)=U(\log H)U^{-1},
\end{equation}
which along with the invariance of the state $\tau^{(2)}$ with
respect to unitary transformations yields
\begin{equation*}
\tau^{(2)}\left[\log\left(UHU^{-1}\right)\right]= \tau^{(2)}[\log H].
\end{equation*}
Hence, \eqref{3.333} implies the equality
\begin{align}\label{1.1.9}
\tau^{(2)}[\log{\bf M}]=
\frac12\left(\tau\left[\log\left(N-KM^{-1}K^*\right)\right]+\tau[\log
M]\right ).
\end{align}

Combining \eqref{3.333} and \eqref{1.1.9} we get
\begin{align*}
&\tau\left[\log\left(M-K^*N^{-1}K\right)\right]+\tau[\log N]\\&=
\tau\left[\log\left(N-KM^{-1}K^*\right)\right]+\tau[\log M].
\end{align*}
By \eqref{Xi:inv} and \eqref{xidef}, this completes the proof of the
theorem.
\end{proof}

\begin{rem}
The requirement that both the Shur complements $M-K^*N^{-1}K$ and
$N-KN^{-1}K^*$ of ${\bf M}$ have bounded inverses is redundant. It is sufficient to
require that at least one of the Shur complements is nonsingular  since
$M-K^*N^{-1}K\in \D_\A$ implies $N-KN^{-1}K^*\in \D_\A$ and vice
versa $($cf. representation \eqref{3.6}$)$.
\end{rem}

The following consequence suggests a recipe  for the computation of the (relative)
Morse index of a $2\times2$ operator matrix (cf. Remark \ref{mors}). It also
provides a representation for the $\xi$-index associated with an
off-diagonal perturbation problem.

\begin{cor}
Assume hypothesis of Theorem \ref{mbs} and let ${\bf M}$ be the operator matrix
$
{\bf M}=\begin{pmatrix}
M& K^*\\
K&N
\end{pmatrix}.
$
 Let $U$ and $W$ be isometries from $\cH$ into $\cH\oplus\cH$ such
that $U^\ast \mathbf{M} U = M$ and $W^\ast \mathbf{M} W = N$. Then
\begin{align}\label{papa}
2\tau^{(2)}[\Xi({\bf
M})]&=\tau\left[\Xi\left (\left(W^*{\bf M^{-1}}W\right)^{-1}\right)\right ]+\tau[\Xi(U^*{\bf M}U)]\\
\nonumber&=\tau\left[\Xi\left (\left(U^*{\bf
M^{-1}}U\right)^{-1}\right)\right ] +\tau[\Xi(W^*{\bf M}W)].
\end{align}
In particular,
\begin{align*}
2\,\xi({\bf M}_0,{\bf M}) & =\xi\left(U^\ast \mathbf{M} U, \left(U^\ast
\mathbf{M}^{-1}U \right)^{-1}\right)\\
& =\xi\left(W^\ast \mathbf{M} W,
\left(W^\ast \mathbf{M}^{-1} W\right)^{-1}\right),
\end{align*}
where ${\bf M}_0=\begin{pmatrix}M& 0\\0& N\end{pmatrix}$.
\end{cor}

\begin{proof}
We notice that equalities \eqref{3.333}  and  \eqref{1.1.9} yield
\begin{align}\label{mama}
2\tau^{(2)}[\Xi({\bf
M})]&=\tau\left[\Xi\left(N-KM^{-1}K^*\right)\right]+\tau[\Xi(M)]\\
\nonumber&=\tau\left[\Xi\left(M-K^*N^{-1}K\right)\right]
+\tau[\Xi(N)].
\end{align}
According to the definition of the
Shur complements of ${\bf M}$ $($cf.  \eqref{shur}$)$,
 we get $$U^\ast \mathbf{M}^{-1} U = (M-K^\ast
N^{-1} K)^{-1},\quad W^\ast \mathbf{M}^{-1} W = (N-K M^{-1}
K^\ast)^{-1}$$ and hence \eqref{papa} follows from \eqref{mama}.
\end{proof}

\begin{rem}\label{xim}
In the $\Iinfty$ setting, a relation similar  to \eqref{papa} has been recently derived
in \cite{NLS}.
\end{rem}

Our next goal is to obtain an extension of the basic invariance
principle stated in Theorem \ref{mbs} by relaxing the invertibility
hypotheses.

\begin{thm}\label{bsp}
Let $M$, $N\in \A$ be dissipative and $K$ an arbitrary operator in
$\A$. Then the following assertions hold.
\begin{enumerate}[(i)]
\item
\begin{align*}\nonumber
&\lim_{\varepsilon\downarrow 0}\xi\left(M+\i\varepsilon
I,M+\i\varepsilon I-K^*(N+\i\varepsilon
I)^{-1}K\right)\\
&=\lim_{\varepsilon\downarrow 0}\xi\left(N+\i\varepsilon
I,N+\i\varepsilon I-K(M+\i\varepsilon I)^{-1}K^*\right).
\end{align*}
\item Assume that $N$ has a bounded inverse. Then
\begin{equation}\label{nach}
\xi\left(M,M-K^*N^{-1}K\right)=\lim_{\epsilon \downarrow
0}\xi\left(N,N-K(M+\ii \epsilon I)^{-1}K^*\right).
\end{equation}
\item If, in addition, the limit
\[K(M+\ii 0 I)^{-1}K^* = \nlim_{\epsilon\downarrow 0} K(M+\i\varepsilon
I)^{-1}K^*\] exists and $N-K(M+\ii 0 I)^{-1}K^*$ has a bounded
inverse, then
\begin{equation}\label{noch}
\xi\left(M,M-K^*N^{-1}K\right)=\xi\left(N,N-K(M+\ii 0
I)^{-1}K^*\right).
\end{equation}
\end{enumerate}
\end{thm}

\begin{proof} (i).
Theorem \ref{mbs} guarantees that
\begin{align*}\nonumber
&\xi\left(M+\i\varepsilon I,M+\i\varepsilon I-K^*(N+\i\varepsilon
I)^{-1}K\right)\\
&\qquad =\xi\left(N+\i\varepsilon I,N+\i\varepsilon
I-K(M+\i\varepsilon I)^{-1}K^*\right),\quad \epsilon>0.
\end{align*}
Therefore, to prove the claim it is sufficient to establish the
existence of the limit
\begin{align*}
&\lim_{\varepsilon\downarrow 0}\xi\left(M+\i\varepsilon
I,M+\i\varepsilon I-K^*(N+\i\varepsilon I)^{-1}K\right)
\\
&=\lim_{\varepsilon\downarrow
0}\bigg(\tau\left[\Xi\left(M+\i\varepsilon I-K^*(N+\i\varepsilon
I)^{-1}K\right)\right]-\tau[\Xi(M+\i\epsilon I)]\bigg).
\end{align*}
By Remark \ref{nepr},  the limit $\lim_{\varepsilon\downarrow
0}\tau[\Xi(M+\i\epsilon I)] $ exists. Next, by Remark \ref{xim},
\begin{equation}\label{kika}
\begin{split}
& \tau\big[\Xi\big(M+\i\varepsilon I-K^*(N+\i\varepsilon
I)^{-1}K\big)\big]\\ & \qquad\qquad =2\tau^{(2)}\big[\Xi({\bf
M}+\i\varepsilon {\bf I})\big]-\tau\big[\Xi(N+\i\varepsilon I)\big],
\end{split}
\end{equation}
where ${\bf M}=\begin{pmatrix}M& K^*\\K& N\end{pmatrix}$ is a
$2\times 2$ operator matrix in $\A\overline{\otimes}M_2$. Applying
Remark \ref{nepr} to the dissipative elements $\bf{M}$ and $N$ in the
algebras $\A\overline{\otimes}M_2$ and $\A$, respectively, insures
the existence of the limit of the left hand side of \eqref{kika} as
$\epsilon\downarrow 0$, completing the proof.

(ii). From Theorem \ref{mbs}, we obtain that the equality
\begin{equation}\label{iixi}
\xi\left(M+\i\epsilon I,M-K^*N^{-1}K+\i\epsilon I\right)=
\xi\left(N,N-K(M+\i\epsilon I)^{-1}K^*\right)
\end{equation}
holds for all $\epsilon>0$. Observe that invertibility of $N$ implies
that of the operator $N-K(M+\i\epsilon I)^{-1}K^*$ for any $\epsilon
> 0$. Indeed, the Herglotz operator-valued function
\[z\mapsto N-K(M+z I)^{-1}K^*\] in the upper-half plane is invertible
for $|\Im z|$ large enough and, therefore, it is invertible for all
$z\in \C_+$ (cf.~\cite[Lemma 2.3]{MakarovP}). Passing to the limit
$\epsilon \downarrow 0$ in \eqref{iixi} and making use of  Remark
\ref{nepr} implies \eqref{nach}.

Since by hypothesis $N-K(M+\ii 0 I)^{-1}K^*$ has a bounded inverse,
using continuity of the operator logarithm (cf. \eqref{intrep}) and that of the state
$\tau$, we attain
\begin{align*}
\lim_{\epsilon \downarrow 0}\tau\left[\Xi\left(N-K(M+\i\epsilon
I)^{-1}K^*\right)\right] =\tau\left[\Xi\left(N-K(M+\ii 0
I)^{-1}K^*\right)\right].
\end{align*}
Now the claim follows from (ii).
\end{proof}

\begin{rem}\label{lemonbound} As one can see from  the proof, for any
dissipative elements $M$ and $N$ in $\A$ and any $K\in \A$, the limit
\begin{equation*}
\lim_{\varepsilon\downarrow 0}\tau\left[\Xi\left(M+\i\epsilon
I-K^*(N+\i\epsilon I)^{-1}K\right)\right]
\end{equation*}
exists. If, in addition, the dissipative operator $M$ has a bounded
inverse, claim (ii) infers the existence of the limit
\begin{equation*}
\lim_{\varepsilon\downarrow 0}\tau\left[\Xi\left(M-K^*(N+\i\epsilon
I)^{-1}K\right)\right].
\end{equation*}
\end{rem}

\section{The Birman-Krein formula revisited}\label{s4}

As an application of Theorem \ref{mbs}, first we state a result
regarding the computation of the relative index associated with
purely imaginary dissipative perturbations $A\mapsto A+\ii B$, $B\geq
0$, of a self-adjoint operator $A$. The following theorem sheds some
light on the role of the characteristic function of a dissipative
operator in the relative index theory. We recall that the Lifshits
characteristic function $\bf S$ of the dissipative operator $A+\ii B$
calculated at the spectral point $\lambda=0$ (see, e.g.,
\cite[Section IV.6]{Gohberg:Krein}) is given by
\begin{equation}\label{S:def}
{\bf S}= I - 2\ii B^{1/2}(A+\ii B)^{-1}B^{1/2}.
\end{equation}

\begin{thm}\label{charfun}
Let  $A=A^*$ and $B=B^*\ge 0$ be elements in $\A$. Suppose both $A$
and $A+\ii B$ have bounded inverses. Then
\begin{align}\label{sss}
\xi(A,A+\ii B)= \frac1\pi \tau\big[\arctan
\big(B^{1/2}A^{-1}B^{1/2}\big)\big]= \frac{1}{2\pi}\tau[\arg\bf S].
\end{align}
Here $\bf S$ is as in \eqref{S:def} and the argument of ${\bf S}$ is
defined by the Spectral Theorem
\begin{equation*}
\arg {\bf S}=\int_{|z|=1} \arg z\,dE_{\bf S}(z),
\quad
\arg \, z\in (-\pi, \pi],\quad z\in \mathbb C\setminus \{0\},
\end{equation*}
with the cut along the negative semi-axis.
\end{thm}

\begin{proof}
Introduce a self-adjoint operator $H=B^{1/2}A^{-1}B^{1/2}$. Theorem
\ref{mbs} implies that
\begin{align}
\xi(A,A+\ii B)&=\xi\left(\ii I,\ii I-H\right) =
\tau\left[ \Xi\left(\ii I-H\right)\right]-\frac{1}{2} \label{1}\\
&=\frac1\pi \tau\left[\Im\log\left(\ii I-H\right)\right]
-\frac{1}{2}\nonumber.
\end{align}
By the Spectral Theorem applied to $H$, we obtain
\begin{align}\label{2}
&\tau\left[\Im\log\left(\ii I-H\right)\right] -\frac{\pi}{2}=
\int_\R\left(\Im\log(\i-\la)-\frac{\pi}{2}\right)
d\tau[E_{H}(\la)]\\\nonumber &= \int_\R\Im\log(1+\i\la)d\tau[E_{
H}(\la)]=\tau\big[\arctan H\big].\nonumber
\end{align}
Combining \eqref{1} and  \eqref{2}, completes the proof of the first
equality in \eqref{sss}.

It is straightforward to verify that
\begin{equation}\label{3}
{\bf S}=(\ii I-H)(\ii I+H)^{-1},
\end{equation}
which, in particular, implies that ${\bf S}$ is unitary. For a smooth
path of unitaries
\begin{equation}\label{pathut}
[0,1]\ni t\mapsto U_t=(\ii I-tH)(\ii I+tH)^{-1}
\end{equation}
linking the identity $I=U_0$ with ${\bf S}=U_1$, we derive that
\begin{align}\label{args}
\tau[\arg {\bf S}]=\Im\,\tau(\log {\bf
S})=\Im\int_0^1\frac{d}{dt}\tau(\log U_t)dt=\Im\int_0^1
\tau\left[\dot U_tU_t^{-1}\right]dt.
\end{align}
Observing that
\begin{equation}\label{1.43}
\begin{split}
\tau\left[\dot U_tU_t^{-1}\right] & =-\tau\big[(H(\ii I+tH)^{-1}\\ &
\qquad\qquad
+(\ii I-tH)(\ii I +tH)^{-2}H)(\ii I+tH)(\ii I-tH)^{-1}\big]\\
&=-\tau\left[H(\ii
I+tH)^{-1}\right]-\tau\left[(\ii I-tH)^{-1}H\right]\\
&=-\frac{d}{dt}\tau[\log(\ii I+tH)]+\frac{d}{dt}\tau[\log(\ii I-tH)],
\end{split}
\end{equation}
we arrive at the equality
\begin{align*}
\Im\int_0^1\tau\left[\dot U_t
U_t^{-1}\right]dt&=\Im\big(\tau[\log(\ii I-H)]-
\tau[\log(\ii I)]\big)\\
&\quad -\Im\big(\tau[\log(\ii I+H)]-\tau[\log(\ii I)]\big).
\end{align*}
By the Spectral Theorem this equality implies
\begin{equation}\label{pri}
\begin{split}
&\Im\int_0^1\tau\left[\dot U_tU_t^{-1}\right]dt\\
&=\int_\R\left(\Im\log\left(\frac{1-\i\la}{\i}\right)-
\Im\log\left(\frac{1+\i\la}{\i}\right)\right)
d\tau[E_H(\la)]\\
&=\int_\R\big(\arg(1+\i\la)-\arg(1-\i\la)\big) d\tau[E_H(\la)]\\
&=2\int\arg(1+\i\lambda)d\tau[E_H(\lambda)]=2\tau[\arctan H].
\end{split}
\end{equation}
Comparing \eqref{args} and \eqref{pri} proves the second equality in
\eqref{sss}.
\end{proof}

Before turning back to the context of perturbation theory for self-adjoint
operators,  it is convenient to collect basic assumptions and related
notation in the form of a hypothesis.
\begin{hyp}\label{final}
Suppose that  $H_0=H_0^*$ and $V=V^*$ are elements in $\A$ and
$H=H_0-V$. Assume  that $V$ is factored in the form $V=-K^*N^{-1}K$,
where $K\in\A$ and $N=N^*$ an element with a bounded inverse in $\A$.
Assume that the norm-limit
\begin{equation*}
K^*(H_0+\ii 0 I)^{-1}K = \nlim_{\varepsilon\downarrow
0}K^*(H_0+\i\varepsilon I)^{-1}K
\end{equation*}
exists and both the operators $\mathcal{N}=N-K^*(H_0+\ii 0 I)^{-1}K$
and $\Re\,\mathcal{N}$ have bounded inverses.

Assume, in addition, that $\bf S$ is the characteristic function of
the dissipative operator $\mathcal{N}$ at the zero value of the
spectral parameter, that is,
\begin{equation}\label{SS:def}
{\bf S}=I - 2\ii (\Im\,\mathcal{N})^{1/2}\mathcal{N}^{-1} (\Im\,
\mathcal{N})^{1/2}.
\end{equation}
\end{hyp}

We conclude (under Hypothesis \ref{final})  with a result relating
the relative index $\xi(H,H_0)$ to the de la Harpe-Skandalis
determinant of the characteristic function $\bf S$  of the
dissipative operator $\mathcal{N}$.

\begin{thm}\label{bkf}
Assume Hypothesis \ref{final}. Let ${\det}_{\tau} {\bf S}$ be the de
la Harpe-Skandalis determinant associated with the homotopy class of
the $C^1$-paths of invertible operators  joining ${\bf S}$ and $I$
and containing the path
\begin{equation}\label{hom}
[0,1]\ni t\mapsto t{\bf S}+(1-t)I.
\end{equation}
Then
\begin{equation}\label{bdfor} {\det}_{\tau} {\bf S}=\Theta\exp \big(-2
\pi \ii \xi(H, H_0)\big),
\end{equation}
where
\begin{equation}\label{teta}
\Theta=\exp \big (-2 \pi \ii
\xi(N, \Re\,\mathcal{N})\big ).
\end{equation}
\end{thm}

\begin{proof}
Applying Theorem \ref{bsp} yields
\begin{equation} \label{bspi1}
\xi(H_0,H)=\xi(N,\mathcal{N})=\xi(N,\Re\,\mathcal{N})+\xi(\Re\,\mathcal{N},\mathcal{N}).
\end{equation}
By Theorem \ref{charfun} one has
\begin{equation*}
\xi(\Re\,\mathcal{N},\mathcal{N})=\frac{1}{2\pi} \tau[\arg {\bf S}]
\end{equation*}
and, hence, \begin{equation} \label{bspi}
\xi(H_0,H)=\xi(N,\Re\,\mathcal{N})+\frac{1}{2\pi }\tau [\arg {\bf S}]
\end{equation}
holds. Multiplying by $2\pi\i$ on both sides of \eqref{bspi} and then
exponentiating ensures, by Lemma \ref{a2} (ii), that
\[\De(t\mapsto U_t)=\Theta\exp\big(-2\pi\ii\xi(H, H_0)\big),\]
where the nonsingular path $[0,1]\ni t\mapsto U_t$ is given by
\begin{align*}
U_t=&\left(\ii
I-t(\im\,\mathcal{N})^{1/2}(\Re\,\mathcal{N})^{-1}(\im\,\mathcal{N})^{1/2}\right)\\
&\times\left(\ii
I+t(\im\,\mathcal{N})^{1/2}(\Re\,\mathcal{N})^{-1}(\im\,\mathcal{N})^{1/2}\right)^{-1},
\quad t\in [0,1].
\end{align*}
Since the path   of unitary operators $ t\mapsto U_t$ with endpoints
${\bf S}$ and $I$ is homotopically equivalent to the path of
invertible operators  $[0,1]\ni t\mapsto t{\bf S}+(1-t)I$  (the point
$-1$ does not belong to the spectrum of ${\bf S}$), the result
follows upon applying Lemma \ref{a1} (i).
\end{proof}

\begin{rem}
Note that the characteristic function $\bf S$ of the dissipative
operator $\mathcal{N}=N-K^*(H_0+\ii 0 I)^{-1}K$ given by
\eqref{SS:def} can also be understood as the abstract scattering
operator associated with the pair $(H_0,H)$ $($cf.~\cite{Yafaev}$)$.
As distinct from the classical Birman-Krein formula \cite{BirmanK}
where the argument of the determinant of the scattering matrix is
directly related to the spectral shift function $($mod $\Z$$)$,
representation \eqref{bdfor} for
 ${\det}_\tau {\bf S}$ via the
$\xi$-index contains a unimodular factor $\Theta$ \eqref{teta}.
Presence  of the additional factor $\Theta$ in \eqref{bdfor} can be
explaiend by the non-integer nature of the $\tau$-Fredholm index for
the pair of orthogonal projections
\begin{equation*}
\xi(N,\Re\,\mathcal{N})=\index_\tau\left(E_N(\R_-),
E_{\Re\,\mathcal{N}}(\R_-)\right)\in [-1,1].
\end{equation*}
\end{rem}

\appendix
\section{}\label{sa}

In this appendix, we recall the concept of a determinant introduced by
P.~de la Harpe and G.~Skandalis in \cite{Skandalis}.

Let $GL^0(\A)$ be the set of boundedly invertible elements of $\A$.
Given a nonsingular $C^1$-path of operators $[0,1]\ni t\mapsto H_t\in
GL^0(\A)$, the de la Harpe-Skandalis determinant associated with the
path $t\mapsto H_t$ is defined by
\begin{equation}\label{uniphase}
\De(t\mapsto H_t)=\exp\left(\int_0^1\tau\left[\dot
H_tH_t^{-1}\right]dt\right).
\end{equation}

Some important properties of the de la Harpe-Skandalis determinant
are listed in the lemma below. The proofs of these facts can be found
in \cite[Lemma 1 and Proposition 2]{Skandalis}.
\begin{lemma}\label{a1}
Suppose that $[0,1]\ni t\mapsto H_t$ is a $C^1$-path of operators in
$GL^0(\A)$.
\begin{enumerate}[(i)]

\item The determinant  $\De(t\mapsto H_t)$ is invariant under fixed
endpoint homotopies.

\item The absolute value of the perturbation determinant
$\De(t\mapsto H_t)$ is path-independent. Moreover,
$$
|\De(t\mapsto H_t)|=\Delta\left(H_1H_0^{-1}\right),
$$
where  $\De(A)=\exp(\tau[\log\sqrt{A^*A}])$ denotes the
Fuglede-Kadison determinant of a boundedly invertible operator
$A\in\A$.

\item If $\|H_t-I\|<1$ for all $t\in [0,1]$, then
\begin{equation}\label{logi}
\De(t\mapsto H_t)=\exp(\tau[\log H_1]-\tau[\log H_0]),
\end{equation}
where the operator logarithm $\log H_j$, $j=1,2,$ in \eqref{logi} is
understood as the norm convergent series
\begin{equation*}
\log H_j=-\sum_{k=1}^\infty \frac{(I-H_j)^k}{k}, \quad j=0,1.
\end{equation*}
\item Let $H^{(j)}_t    :[0,1]\to GL^0(\A)$, $j=1,2$, be $C^1$-paths.
Then
\begin{equation*}
\De\left(t\mapsto H^{(1)}_tH^{(2)}_t\right)=\De\left(t\mapsto
H_t^{(1)}\right) \De\left(t\mapsto H_t^{(2)}\right).
\end{equation*}
\end{enumerate}
\end{lemma}

The following result reduces the computation of the determinant for
paths of operators in either $\D_\A$ or
\begin{equation*} \mathcal{U}_\A = \{U\,
:\, U=(\i I-H)(\i I+H)^{-1}\quad \text{for some}\quad H=H^*\in\A\}
\end{equation*}
to that of the state $\tau$ of the operator logarithm.

\begin{lemma}\label{a2}
(i) For a $C^1$-path of operators $[0,1]\ni t\mapsto H_t\in\D_\A$
with $H_0=I$,
\begin{equation}\label{detn}
\De(t\mapsto H_t)=\exp(\tau[\log H_1]),
\end{equation}
where $\log(\cdot)$ is the principal branch of the operator logarithm
of $H_1$ with the cut along the negative imaginary semi-axis provided
by the Riesz functional calculus.

(ii) For a $C^1$-path of operators $[0,1]\ni t\mapsto
U_t\in\mathcal{U}_\A$ with $U_0=I$,
\begin{equation*}
\De(t\mapsto U_t)=\exp\left(\tau\left[\widetilde{\log}\,
U_1\right]\right),
\end{equation*}
where $\widetilde{\log}(\cdot)$ is the principal branch of the
operator logarithm of $U_1$ with the cut along the negative real
semi-axis provided by the Spectral Theorem.
\end{lemma}

\begin{proof} (i) One notices that $\tau\left[\dot
H_tH_t^{-1}\right]=\frac{d}{dt}\tau[\log H_t]$. Integrating the
latter expression from $0$  to $1$ and comparing the result with
\eqref{uniphase} implies \eqref{detn}. The proof of (ii) goes along
the same lines as that of (i).
\end{proof}


\end{document}